\documentclass[a4paper,12pt]{article}

\usepackage{url, color, epsfig, amsmath, amsthm, amssymb}

\usepackage{hyperref}

\usepackage[margin = 1.0in]{geometry}

\usepackage{algorithm}
\usepackage[noend]{algorithmic}
\usepackage{appendix}
\allowdisplaybreaks

\newtheorem{theorem}{Theorem}[section]
\newtheorem{corollary}[theorem]{Corollary}
\newtheorem{lemma}[theorem]{Lemma}

\newtheorem{proposition}[theorem]{Proposition}

\renewcommand{\Re}{\mathbb{R}}
\newcommand{\Ex}{\mathbb{E}}

\newcommand{\F}{\mathcal{F}}

\renewcommand{\H}{\mathcal{H}}

 \newcommand{\G}{\mathcal{G}}

\newcommand{\N}{\mathbf N}

\renewcommand{\S}{\mathcal S}

\newcommand{\al}{\alpha}
\newcommand{\be}{\beta}
\newcommand{\ga}{\gamma}
\newcommand{\tri}{\triangle}
\newcommand{\eps}{\epsilon}

\def\T{\ensuremath{\mathcal{T}}}

\newcommand{\remove}[1]{}

\renewcommand{\S}{\mathcal{S}}

\DeclareMathOperator{\Ave}{\texttt{Ave}}
\DeclareMathOperator{\conv}{\texttt{conv}}

\DeclareMathOperator{\pos}{\texttt{pos}}
\DeclareMathOperator{\diam}{\texttt{diam}}
\DeclareMathOperator{\centroid}{\texttt{c}}

\def\db2{\ensuremath{\lfloor d/2 \rfloor}}

\newcommand{\ptas}[1]{\textsc{PTAS}}

\def\F{\ensuremath{\mathcal{F}}}
\def\H{\ensuremath{\mathcal{H}}}

\def\eps{\ensuremath{\epsilon}}
\def\de{\ensuremath{\delta}}

\renewcommand{\Re}{\mathbb{R}}

\title{Theorems of Carath\'eodory, Helly, and Tverberg without dimension}

\author{Karim Adiprasito, Imre B\'ar\'any, Nabil H. Mustafa, and Tam\'as Terpai}

\date{}

\begin{document}

\maketitle

\begin{abstract}  We initiate the study of no-dimensional versions of classical theorems in convexity. One example is Carath\'eodory's theorem without dimension: given an $n$-element set $P$ in a Euclidean space, a point $a \in \conv P$, and an integer $r \le n$, there is a subset $Q\subset P$ of $r$ elements such that the distance between $a$ and $\conv Q$ is less than $\diam P/\sqrt {2r}$. The similar no-dimension Helly theorem states that given $k\le d$ and a finite family $\F$ of convex bodies, all contained in the Euclidean unit ball of $\Re^d$, there is a point $q \in \Re^d$ which is closer than $1/\sqrt k$ to every set in $\F$. This result has several colourful and fractional consequences. Similar versions of Tverberg's theorem and some of their extensions are also established.
\end{abstract}

\section{Introduction}

Two classical results, theorems of Carath\'eodory and Helly, both more than hundred years old, lie at the heart of combinatorial convexity. Another basic result is Tverberg's theorem (generalizing Radon's), it is more than fifty years old and is equally significant. In these theorems dimension plays an important role. What we initiate in this paper is the study of the dimension free versions of these theorems. 

Let us begin with Carath\'eodory's theorem \cite{Carath}. It says that every point in the convex hull of a point set $P \subset \Re^d$ is in the convex hull of a subset $Q\subset P$ with at most $d+1$ points. Can one require here that $|Q|\le r$ for some fixed $r\le d$? The answer is obviously no. For instance when $P$ is finite, the union of the convex hull of all $r$-element subsets of $P$ has measure zero while $\conv P$ may have positive measure. So one should set a more modest target. One way for this is to try to find, given $a \in \conv P$, a subset $Q \subset P$ with $|Q|\le r$ so that $a$ is close to $\conv Q$. This is the content of the following theorem:

\begin{theorem} \label{th:Cara}
Let $P$ be a set of $n$ points in $\Re^d$, $r \in [n]$ and $a \in \conv P$. Then there
exists a subset $Q$ of $P$ with $|Q|=r$ such that
\begin{equation}
d\left( a, \conv Q \right) < \frac {\diam P}{\sqrt {2r}}.
\end{equation}
\end{theorem}

When $r\ge d+1$, the stronger conclusion $a \in \conv Q$ follows of course from Carath\'eodory's theorem. But in the statement of the theorem the dimension $d$ has disappeared. So one can think of the $n$-element point set $P$ as a set in $\Re^n$ (or $\Re^{n-1}$) with $a \in \conv P$. The conclusion is that for every $r<n$ the set $P$ has a subset $Q$ of size $r$ whose convex hull is close to $a$. That is why we like to call the result ``no-dimension Carath\'eodory theorem''. 

Next comes Helly's theorem. It is about a family of convex sets $K_1,\ldots,K_n$ in $\Re^d$, $n\ge d+1$. For $J \subset [n]$ define $K(J)=\bigcap_{j\in J}K_j$. With this notation Helly's theorem says that if $K(J)\ne \emptyset$ for every $J\subset [n]$ with $|J|=d+1$, then $\bigcap_1^nK_j\ne \emptyset$. Again, what happens if the condition $K(J)\ne \emptyset$ only holds when $|J|=k$ and $k \le d$? Then the statement fails to hold for instance when each $K_i$ is a hyperplane (and they are in general position). But again, something can be saved, namely:

\begin{theorem}\label{th:helly} Assume $K_1,\ldots,K_n$ are convex sets in $\Re^d$ and $k \in [n]$. For $J \subset [n]$ define $K(J)=\bigcap_{j \in J}K_j$. If the Euclidean unit ball $B(b,1)$ centered at $b \in \Re^d$ intersects $K(J)$ for every $J\subset [n]$ with $|J|=k$, then  there is point $q \in \Re^d$ such that
\[d(q,K_i) < \frac 1 {\sqrt k} \mbox{ for all }i \in [n].
\]
\end{theorem}
Again, the dimension of the underlying space is not important here. Note however that Helly's theorem is invariant under (non-degenerate) affine transformations while its dimension free version is not. The same applies to the no-dimension Carath\'eodory theorem. We mention that the condition $K_i \subset B(b,1)$ for all $i\in [n]$ could be replaced by $K_i \subset B(b,r)$, and then the statement would be $d(q,K_i) < \frac r {\sqrt k}$ for all $i \in [n]$.

The dimension free version of Tverberg's famous theorem~\cite{tver} is as follows.
\begin{theorem}
\label{th:nodimensiontverberg}
Given a set $P$ of $n$ points in
$\Re^d$ and an integer $2 \le k \le n$, there exists a point $q \in \Re^d$ and a partition of $P$ into $k$ sets  $P_1, \ldots, P_k$ such that
\[
d \left(q, \conv P_i \right) \leq (2+\sqrt 2) \cdot \sqrt{\frac kn} \diam P \qquad \mbox{ for every } i\in [k].
\]
\end{theorem}

The goal of this paper is to prove these results and their more general versions together with several applications. We expect them to be highly useful just as their classical versions have been.

\section{Carath\'eodory without dimension}

In Theorem~\ref{th:Cara} the appearance of the scaling factor $\diam P$ is quite natural. The dependence on $r$ is best possible: when $d=n-1$ and $P$ is the set of vertices of a regular $(n-1)$-dimensional simplex whose centre is $a$, then for every $Q\subset P$ with $|Q|=r$,
\[
d(a,\conv Q)= \sqrt{\frac 1{2r}-\frac 1{2n}}\diam P,
\]
which is asymptotically the same as the upper bound in Theorem~\ref{th:Cara} in the no dimension setting.

The coloured version of Carath\'eodory's theorem \cite{bar} states that if $a \in \bigcap_1^{d+1} \conv P_i$, where $P_i\subset \Re^d$, then there is a transversal $T=\{p_1,\ldots,p_{d+1}\}$ such that $a \in \conv T$. Here a {\sl transversal} of the set system $P_1, \ldots, P_{d+1}$ is a set $T=\{p_1,\ldots,p_{d+1}\}$ such that $p_i \in P_i$ for all $i \in [d+1]$. We extend this to the no-dimension case as follows.

\begin{theorem}\label{th:CalCar}
Let $P_1, \ldots, P_r$ be $r\ge 2$ point sets in $\Re^d$ such that $a \in \bigcap_1^r \conv P_i$. Define $D = \max_{i\in [r]} \diam P_i$. Then there exists a transversal $T$ such that
$$d\left( a, \conv T \right)<  \frac D{\sqrt{2r}}.$$
\end{theorem}


We remark further that this result, in a somewhat disguised form, has been know for quite some time. This is explained in the next section. 

The proof is an averaging argument that can be turned into a randomized algorithm that finds the transversal  $T$ in question; the method of conditional probabilities also gives a deterministic algorithm.  We mention that a recent paper of Barman~\cite{barman} proves a qualitatively and quantitatively weaker statement, applicable only
for the case $r = d+1$: given $d+1$ point sets $P_1, \ldots, P_{d+1}$ with  $a \in \bigcap_1^{d+1} \conv P_i$, it is shown how to compute, using
convex programming, a subset of $r$ points $P'$ with $|P' \cap P_i| \leq 1$ for each $i \in [d+1]$, such
that $d\left(a, \conv P'\right) = O\left( \frac{D}{\sqrt{r}}\right)$.
We improve on this in two ways: $a)$ the parameter $r$, the number of sets $P_i$, can be any value $r \leq d$, and thus truly does not depend on the dimension,
and $b)$ the running time of Barman's algorithm is $(nd)^{O(r)}$ while the one from Theorems~\ref{th:CalCar} and \ref{th:CalCar2} is $O(nd)$.
We remark further that in the case $r = d+1$, finding the transversal $T$ such that $a \in \conv T$ in time polynomial in the number $n$ of the input points
and the dimension is a longstanding open problem (see~\cite{Mat}).
 Barman's work implies an algorithm that computes an approximate transversal
with running time $(nd)^{O(r)}$, while Theorem~\ref{th:CalCar} improves the running time to $O(nd)$.

A strengthening  of the colourful Carath\'eodory's theorem from \cite{ABBFM} and \cite{HoPaTv} states that given non-empty sets $P_1,\ldots, P_{d+1} \in \Re^d$ such that $a \in \conv (P_i \cup P_j)$ for every $i,j \in [d+1], i\ne j$, there is a transversal $T=\{p_1,\ldots,p_{d+1}\}$ such that $a \in \conv T$. It is shown in \cite{ABBFM} that the ``union of any two'' condition here cannot be replaced by the ``union of any three'' (or more) condition. We extend this result to the no-dimensional case with ``the union of any two or more" condition.

\begin{theorem} \label{th:CalCar2}
Let $P_1, \ldots, P_r$ be $r\ge 2$ point sets in $\Re^d$, $D = \max_{i\in [r]} \diam P_i$, and $t\in [r-1]$. Assume that
for distinct $i_1,i_2,\ldots,i_t \in [r]$ we have $a \in \conv (P_{i_1}\cup\ldots \cup P_{i_t})$. Then there exists a transversal $T=\{p_1,\ldots,p_r\}$ such that
$$d\left( a, \conv T \right) \leq \be \cdot \frac{\diam P}{\sqrt{r-t+1}},$$
where $\be = 4 \sqrt{ \frac{ \ln 4}{3}} =2.71911...$.
\end{theorem}

The proof is based on the Frank-Wolfe procedure~\cite{bhpr, clark, FrWo}. For the case $t =1 $, it implies a slightly
weaker bound than Theorem~\ref{th:CalCar}, i.e. it finds a transversal $T$ with
\begin{equation}\label{eq:detalg}
d\left( a, \conv T \right) \leq \be \cdot \frac{D}{\sqrt{r}}.
\end{equation}

\medskip
There is a cone version of Carath\'eodory's theorem which is stronger than the convex version. Writing $\pos P$ for the cone hull of $P\subset \Re^d$, it says the following. Assume $P \subset \Re^d$ and $a\in \pos P$ and $a\ne o$. Then there is $Q \subset P$ with $|Q|\le d$ such that $a \in \pos Q$. The corresponding no-dimension variant would say that under the same condition and given $r<d$, there is $Q \subset P$ with $|Q|\le r$ such that the angle between $a$ and the cone $\pos Q$ is smaller than some function of $r$ that goes to zero as $r\to \infty$. Unfortunately, this is not true as the following example shows.

{\bf Example.} Let $P=\{v_1,\ldots,v_d\}$ be the set of vertices of a regular $(d-1)$-dimensional simplex. Assume that its centre of gravity, $a$, is the closest point of $\conv P$ to the origin, and $|a|=h$ is small.
Then $a \in \pos P$. For any subset $Q$ of $P$, of size $r<d$, $\pos Q$ is contained in the boundary of the cone $\pos P$. The minimal angle $\phi$ between $a$ and a vector on the boundary of $\pos Q$ satisfies
\[\
\tan \phi = \frac {|v_1-a|}{(d-1)h}
\]
and can be made arbitrarily large by choosing $h$ small enough.

\section{Earlier and related results}\label{sec:hist}
Results similar to Theorem~\ref{th:Cara} but without the no-dimension philosophy have been known for some time. Each comes with a different motivation. The first seems to be the one by Starr~\cite{Sta} from 1969, see also ~\cite{Sta2}. It measures the non-convexity of the set $P_1+\ldots+P_r$ and is motivated by applications in economy. The result is almost the same as Theorem~\ref{th:CalCar}, only the scaling factor is different. The proof uses the Folkman-Shapley lemma (c.f.~\cite{Sta}). A short and elegant proof using probability is due to Cassels~\cite{Cas}. A comprehensive survey of this type results is given in~\cite{FMMZ}. 

A similar but more general theorem of B. Maurey appeared in 1981 in a paper of Pisier~\cite{pisi}. It is motivated by various questions concerning  Banach spaces. It says that if a set $P$ lies in the unit ball of the space, $a \in \conv P$ and $r\in \N$, then $a$ is contained in a ball of radius $\frac  c{\sqrt r}$ whose center is the centroid of a multiset $Q\subset P$ with exactly $r$ elements, where $c$ is a constant. Here $P$ is supposed to be $0$-symmetric (or $\conv P$ is the ``absolute convex hull''). The proof uses Khintchin's inequality and is probabilistic. Further results of this type were proved by Carl~\cite{Carl} and by Carl and Pajor~\cite{CaPa} and used in geometric Banach space theory. Unlike in Starr's theorem, the underlying space is not necessarily Euclidean, for instance $L_p$ spaces are allowed. Some of the results in this area have become highly influential in geometric concentration of measure (see \cite{Gue, GiMi} for an overview) as for instance Talagrand's inequality (convex subsets of the cube of some measure are highly exhaustive), which is dimension independent as well.

Another way of stating Theorem~\ref{th:Cara} is this. It is possible to find, given a parameter $\eps > 0$, $O\left(\frac{1}{\eps^2}\right)$ points of $P$ whose convex hull is within distance $\eps \cdot \diam P$ from $a$. Such a result was discovered in 2015 by S. Barman~\cite{barman}. His proof is almost identical to that of Maurey or Pisier~\cite{pisi}. But the motivation there is very different. In fact Barman~\cite{barman} has found a beautiful connection of such a statement to additive approximation algorithms. Similar connection was established in ~\cite{AlonLee} as well. The basic idea is the following. Consider an optimization problem that can be written as a bilinear program---namely maximizing/minimizing an objective function of the form $x^T A y$, where the variables are $x, y \in \Re^n$.
If one knew the optimal value of the vector $y$, then the above bilinear program reduces to a linear one, which can be solved in polynomial time. Barman showed that several problems---among them computing Nash equilibria and densest bipartite subgraph problem---have two additional properties: $i)$ $y$ lies inside the convex-hull of some polytope, and $ii)$ if $y$ and $y'$ are two close points in $\Re^n$, then the value of the bilinear programs on $y$ and $y'$ are also close. Then applying the above approximate version of Carath\'eodory's theorem for the optimal point $y$ (whose actual value we don't know), there must exist a point $y'$, depending on a $O(\frac{1}{\eps^2})$-sized subset of the input, such that the distance between $y$ and $y'$ is small. Now one can enumerate all $O(\frac{1}{\eps^2})$-sized subsets to compute all such $y'$, and thus arrive at an approximation to the bilinear program.

A similar inequality was proved by B\'ar\'any and F\"uredi~\cite{BF} in 1987 with a very different purpose. They showed that every deterministic polynomial time algorithm that wants to compute the volume of a convex body in $\Re^d$ has to make a huge error, namely, a multiplicative error of order $\left(\frac d{\log d}\right)^d$. Their proof is based on a lemma similar to Theorem~\ref{th:Cara}. Before stating it we have to explain what the $\rho$-{\sl cylinder} above a set $Q\subset \Re^d$ is, where $|Q|\le d$. Let $B$ denote the Euclidean unit ball of $\Re^d$, and let $L$ be the linear (complementary)  subspace orthogonal to the affine hull of $Q$. Then the cylinder in question is $Q^{\rho}:=(L\cap \rho B)+\conv Q$. With this notation the key lemma in \cite{BF} says that given $P\subset B$ and $r\le d$, every point in $\conv P$ is contained in a cylinder $Q^{\rho(d,r)}$ for some $Q\subset P$ of size $r$, here $\rho(d,r)=\sqrt {\frac {d-r+1}{d(r-1)}}$.  This becomes $\frac{1}{\sqrt {r-1}}$ in the no-dimension setting as
\[
\sqrt {\frac {d-r+1}{d(r-1)}}=\sqrt {\frac 1{r-1}-\frac 1d}<\sqrt {\frac 1{r-1}}.
\]
 and would give in Theorem~\ref{th:Cara} the estimate
\[
d\left( a, \conv Q \right) \leq  \frac {R}{\sqrt {r-1}}
\]
where $R$ is the radius of the ball circumscribed to $P$. By Jung's theorem~\cite{jun}, $R\le \sqrt {\frac {n-1}{2n}}D$, which gives in our setting the slightly weaker upper bound
\[
d\left( a, \conv Q \right) \leq  \frac {D}{\sqrt {2(r-1)}}.
\]
The proof of the lemma from \cite{BF} does not seem to extend to the case of Theorem~\ref{th:CalCar}.

We note that the estimates in Starr's theorem, in Maurey's lemma (and Barman's), and the one in \cite{BF}, and also in Theorems~\ref{th:Cara} and ~\ref{th:CalCar} are all of order $\frac{1}{\sqrt r}$ but the constants are different. Part of the reason is that the setting is slightly different: in the first ones $P$ is a subset of the unit ball of the space while in Theorem~\ref{th:Cara} and ~\ref{th:CalCar} (and elsewhere in this paper) the scaling parameter is $\diam P$.

\section{Helly's theorem without dimension}
\label{sec:helly}

The no-dimension version of Helly, Theorem~\ref{th:helly} extends to the colourful version of Helly's theorem, which is due to Lov\'asz and which appeared in \cite{bar}, and to the fractional Helly theorem of Katchalski and Liu~\cite{KaLi}, cf~\cite{Kalai}. Their proofs are based on a more general result. To state it some preparation is needed. We let $B$ or $B^d$ denote the (closed) Euclidean unit ball in $\Re^d$ and write $B(a,\rho)$ for the Euclidean ball centred at $a\in \Re^d$ of radius $\rho$. Suppose $\F_1,\ldots,\F_k$ are finite and non-empty families of convex sets in $\Re^d$, $\F_i$ can be thought of as a collection of convex sets  of colour $i$. A {\sl transversal} $\T$ of the system $\F_1,\ldots,\F_k$ is just $\T=\{K_1,\ldots,K_k\}$ where $K_i \in \F_i$ for all $i \in [k]$. We define $K(\T)=\bigcap_1^kK_i$. Given $\rho_i > 0$ for all $i \in [k]$, set $\rho=\sqrt{\rho_1^2+\ldots +\rho_k^2}$.

\begin{theorem}\label{th:genhelly} Assume that, under the above conditions, for every $p \in \Re^d$ there are at least $m_i$ sets $K \in \F_i$ with $B(p,\rho_i)\bigcap K=\emptyset$ for all $i \in [k]$. Then for every $q \in \Re^d$ there are at least $\prod_1^k m_i$ transversals $\T$ such that
\[
d(q,K(\T)) >\rho,
\]
with the convention that $d(q,\emptyset)=\infty$.
\end{theorem}

We mention that the value $\prod_1^km_i$ is best possible as shown by the following example. Let $e_1,\ldots,e_k$ denote the standard basis vectors of $\Re^k$ and choose a real number $r_i$ larger than $\rho_i$, but only slightly larger. Set $v_i=r_ie_i$. For every $i \in [k]$ the family $\F_i$ contains $m_i$ copies of the hyperplane $H_i^-=\{x \in \Re^k: v_i(x-v_i)=0\}$ and also $m_i$ copies of the hyperplane $H_i^+=\{x \in \Re^k: v_i(x+v_i)=0\}$, and furthermore some finitely many copies of the whole space $\Re^k$. It is clear that the smallest ball intersecting every set in $\F_i$ is $r_iB^k$. Moreover, given a transversal $H_1^{\eps_1},\ldots,H_k^{\eps_k}$ of the system $\F_1,\ldots,\F_k$ with $\eps_i \in \{+,-\}$, their intersection is a point at distance $r=\sqrt{r_1^2+\ldots +r_k^2}$ from the origin, and there are exactly $\prod_1^km_i$ such transversals. All other transversals of the system have a point in the interior of $r B$, and then also in $\rho B$ if the $r_i$s are chosen close enough to $\rho_i$.

Here comes the no-dimension colourful variant of Helly's theorem.

\begin{theorem}\label{th:Colhelly} Let $\F_1,\ldots,\F_k$ be finite and non-empty families of convex sets in $\Re^d$. If for every transversal $\T$ the set $K(\T)$ intersects the Euclidean unit ball $B(b,1)$, then there is $i\in [k]$ and a point $q \in \Re^d$ such that
\[d(q,K)\le \frac 1 {\sqrt k} \mbox{ for all }K \in \F_i.
\]
\end{theorem}

The proof is just an application of Theorem~\ref{th:genhelly} with $\rho_i=\frac 1{\sqrt k}$: if for every $q \in \Re^d$ and for every $i\in [k]$ there is a $K\in \F_i$ with $d(q,K)>1/\sqrt k$, then $m_i\ge 1$ for all $i \in [k]$. And the theorem implies the existence of a transversal with $d(q,K(\T)) >\rho$. For the fractional version set $|\F_i|=n_i$.

\begin{theorem}\label{th:Colfrachelly} Let $\al \in (0,1]$ and define $\be=1-(1-\alpha)^{\frac{1}{k}}$. Assume that for an $\al$ fraction of the transversals $\T$ of the system $\F_1,\ldots,\F_k$, the set $K(\T)$ has a point in $B(b,1)$. Then there is $q \in \Re^d$ and $i\in [k]$ such that at least $\be n_i$ elements of $\F_i$ intersect the ball $B(q,1/\sqrt k)$.
\end{theorem}

An example similar to the above one shows that the value $\be=1-(1-\alpha)^{\frac{1}{k}}$ is best possible.

Theorem~\ref{th:Colfrachelly} is a consequence of Theorem~\ref{th:genhelly} again. Indeed, if no ball $B(q,1/\sqrt k)$ intersects $\be n_i$ elements of $\F_i$, then $m_i>(1-\be)n_i$. If this holds for all $i\in [k]$, then the number of transversals $\T$ that are disjoint from a fixed unit ball is larger than
$$(1-\be)^k\prod_1^kn_i < \prod_1^km_i \le (1-\al)\prod_1^kn_i$$
contrary to the assumption that $K(\T)$ intersect $B(b,1)$ for an $\al$ fraction of the transversals.

The case when all $\F_i$ coincide with a fixed family $\F=\{K_1,\ldots,K_n\}$ is also interesting and a little different because the transversals correspond to $k$-tuples from $\F$ with possible repetitions. But the proof of Theorem~\ref{th:genhelly} can be modified to give the following result.

\begin{theorem}\label{th:frachelly} Again let $\al \in (0,1]$ and define $\be=1-(1-\alpha)^{\frac{1}{k}}$. Let $\F$ be a finite family of convex sets in $\Re^d$, $|\F|\ge k$. Assume that for an $\al$ fraction of $k$-tuples $K_1,\ldots,K_k$ of $\F$, the set $\bigcap_1^kK_i$ has a point in $B(b,1)$. Then there is $q \in \Re^d$ such that at least $\be |\F|$ elements of $\F$ intersect the ball
$B(q,1/\sqrt k)$.
\end{theorem}

\section{Further results around Helly's theorem}
\label{sec:hellymore}

A more precise version of Theorem~\ref{th:helly} is the following one.
\begin{theorem}\label{th:hellyprec}
Under the conditions of Theorem~\ref{th:helly} there is point $q \in \Re^d$ such that
\begin{equation}\label{eq:simplex}
d(q,K_i)\le \sqrt{\frac {n-k}{k(n-1)}} \mbox{ for all } i \in [n].
\end{equation}
\end{theorem}

This bound is best possible, as shown by a regular simplex $\tri$ on $n$ vertices whose inscribed ball is $B(b,r)$ where $r=\sqrt{\frac {n-k}{k(n-1)}}$: let $K_i$ be the closed halfspace such that $K_i\cap \tri$ is the $i$th facet of $\tri$ ($i \in [n]$) and set $\F=\{K_1,\ldots,K_n\}$. Direct computation shows then that the ball $B(b,1)$ has a single point in common with $K(J)$ for every $J\subset [n]$. This example also shows that the bound in Theorem~\ref{th:helly} is best possible in the no-dimension setting as
\[
 \sqrt{\frac {n-k}{k(n-1)}}=\sqrt{\frac 1k-\frac 1{n-1}+\frac 1{k(n-1)}} = \sqrt{ \frac 1k} +o(k).
\]

The proof is based on a geometric inequality about simplices. It says the following.

\begin{theorem}\label{th:simplex} Let $\tri$ be a (non-degenerate) simplex on $n$ vertices with inradius $r$ and let $k\in [n]$. Then any ball intersecting the affine span of each $(k-1)$-dimensional face of $\tri$ has radius at least $\lambda_n r$ where $\lambda_n=\sqrt{\frac{(n-1)(n-k)}{k}}$ is the optimal ratio for the regular simplex.
\end{theorem}
The case $k=n-1$ is a tautology, and the case $k=1$ is well-known: it is just the fact that the radius of the circumscribed ball is at least dimension times the inradius. To our surprise we could not find the general case in the literature, even in the weaker form saying that any ball intersecting each $(k-1)$-dimensional face of $\tri$ has to have radius at least $\lambda_n r$. Theorems~\ref{th:hellyprec} and \ref{th:simplex} are not directly connected to our no-dimensional setting. That's why their proofs will to be given separately, in the last section.

Theorem~\ref{th:genhelly} becomes completely trivial when any $\rho_i=0$. But stronger statements hold when $k=d+1$ and every $\rho_i=0$, namely, one can show that $K(\T)=\emptyset$ for a certain number of transversals.  For instance, when each $m_i=1$, which simply means that no family $\F_i$ is intersecting, we can show that there is a transversal with $K(\T)=\emptyset$. This is exactly the colourful Helly theorem of Lov\'asz. When $\F_1=\ldots=\F_{d+1}$, then the corresponding statement is just Helly's theorem. The proofs of Theorems~\ref{th:genhelly} and \ref{th:frachelly} can easily be modified to give (another) new proof the (colourful) Helly theorem.

Perhaps the most interesting case is the fractional and colourful Helly theorem saying that, for every $\al\in (0,1]$ and every $d$ there is $\be=\be(\al,d)>0$ such that the following holds. If $d+1$ finite families $\F_1,\ldots,\F_{d+1}$ of convex sets in $\Re^d$ are given so that an $\al$ fraction of their transversals intersect, then $\F_i$ contains an intersecting subfamily of size $\be|\F_i|$ for some $i \in [d+1]$. (The no-dimension version is Theorem~\ref{th:Colfrachelly} above.) This result was first stated and proved in ~\cite{BFMOP} with $\be=\al/(d+1)$.
The better bound $\be = \max \{ \al/(d+1),1-(d+1)(1-\al)^{1/(d+1)}\}$ was obtained by M Kim~\cite{MKim} and is the best lower bound at present.

The following example shows that $\be \le 1-(1-\al)^{1/(d+1)}$. $\F_i$ consists of $m_i=(1-\be)n_i$ parallel hyperplanes (in $\Re^d$) and $n_i-m_i$ copies of  $\Re^d$ (for every $i$). The hyperplanes are chosen so that if $\T$ is a transversal of hyperplanes only, then $K(\T)=\emptyset$. Then the number of non-intersecting transversals is exactly \[\prod_1^{d+1}m_i=(1-\be)^{d+1}\prod_1^{d+1}n_i=(1-\al)\prod_1^{d+1}n_i,
\] and the largest intersecting subfamily of $\F_i$ is of size $n_i-m_i+1=\be n_i+1$. So the question is whether this upper bound is tight. We hope to return to this in a subsequent paper.

We remark that in the original fractional Helly theorem of Katchalski and Liu~\cite{KaLi} a single family $\F$ of convex sets is given with the property that an $\al$ fraction of the $d+1$-tuples of $\F$ are intersecting. And the conclusion is that $\F$ contains an intersecting subfamily of size $\be|\F|$ with $\be=\al/(d+1)$. A famous result of Kalai~\cite{Kalai} from 1984 shows that $\be\ge 1-(1-\al)^{1/(d+1)}$, which is best possible.

\medskip

\section{Tverberg's theorem without dimension}
\label{sec:tverberg}

We are going to prove the no-dimension version of the more general coloured Tverberg Theorem (cf. \cite{ZiVr} and \cite{BlaZie}). We assume that the sets $C_1,\ldots,C_r\subset \Re^d$ (considered as colours) are disjoint and each has size $k$. Set $P=\bigcup_1^rC_j$.

\begin{theorem}\label{th:nodimColTv} Under the above conditions there is a point $q \in \Re^d$ and a partition $P_1,\ldots,P_k$ of $P$ such that $|P_i\cap C_j|= 1$ for every $i\in[k]$ and every $j \in [r]$ satisfying
\[
d(q,\conv P_i)\le (1+\sqrt 2)\frac {\diam P}{\sqrt r} \mbox{ for every } i \in [k].
\]
\end{theorem}

This result implies the uncoloured version, that is, Theorem~\ref{th:nodimensiontverberg}. To see this we write $|P|=n=kr+s$ with $k\in \N$ so that $0\le s\le r-1$. Then delete $s$ elements from $P$ and split  the remaining set into sets (colours) $C_1,\ldots,C_r$, each of size $k$. Apply the coloured version and add back the deleted elements (anywhere you like). The outcome is the required partition, the extra factor $\sqrt 2$ between the constants $2+\sqrt 2$ and $1+\sqrt 2$ comes when $k=2$ and $2r$ is only slightly smaller than $n$. But Theorem~\ref{th:nodimensiontverberg} holds with constant $1+\sqrt 2$ (instead of $2+\sqrt 2$) when $r$ divides $n$.

We remark further that the bounds given in Theorems~\ref{th:nodimensiontverberg} and ~\ref{th:nodimColTv} are best possible apart from the constants. Indeed, the regular simplex with $n=kr$ vertices shows, after a fairly direct computation, that for every point $q \in \Re^{n-1}$ and every partition $P_1,\ldots,P_k$ of the vertices
\[
\max_{i \in [k]} d(q,\conv P_i)\ge \frac {\diam P}{\sqrt {2r}} \sqrt{\frac {k-1}k\frac {n-1}{n-2}}.
\]
The computation is simpler in the coloured case. We omit the details.

\section{Applications}\label{sec:appl}

Several applications of the Carath\'edory, Helly, and Tverberg theorems extend to the no-dimension case. We do not intend to list them all.
But here is an example: the centre point theorem of Rado~\cite{Rad} saying that given a set $P$ of $n$ points in $\Re^d$,
there is a point $q \in \Re^d$ such that any half-space containing $q$ contains at least $\frac{n}{d+1}$ points
of $P$. The proportion $\frac{1}{d+1}$ cannot be improved, in the sense that there exist examples where every point in $\Re^d$
has some half-space containing it and containing at most $\lceil \frac{n}{d+1} \rceil$ points of $P$. The no-dimension version
goes beyond this---at the cost of approximate inclusion by half-spaces.

\begin{theorem}
\label{th:centerpoints}
Let $P$ be a set of $n$ points in $\Re^d$  lying in the unit ball $B(b, 1)$. For any integer $k > 0$, there exists a point $q \in \Re^d$  such that
any closed half-space containing $B\left(q, \frac{1}{\sqrt{k}} \right)$ contains at least $\frac{n}{k}$ points of $P$.
\end{theorem}

The proof is easy and is omitted.

We also give no-dimension versions of the selection lemma \cite{bar} and \cite{Mat} and the weak $\eps$-net theorem \cite{ABFK} .

\begin{theorem}\label{th:selection} Given a set $P\subset \Re^d$ with $|P|=n$ and $D=\diam P$ and an integer $r \in [n]$, there is a point $q \in \Re^d$ such that the ball $B\left(q,\frac {3.5D}{\sqrt r}\right)$ intersects the convex hull of $r^{-r} {n \choose r}$ $r$-tuples in $P$.
\end{theorem}

As expected, the no-dimension selection lemma implies the weak $\eps$-net theorem, no-dimensional version.

\begin{theorem}\label{th:epsnet} Assume $P \subset \Re^d$, $|P|=n$, $D=\diam P$, $r\in [n]$ and $\eps >0$. Then there is a set $F \subset \Re^d$ of size at most $r^r \eps ^{-r}$ such that for every $Y\subset P$ with $|Y| \ge \eps n$
\[
\left(F+\frac {3.5D}{\sqrt r}B\right)\cap \conv Y \ne \emptyset.
\]
\end{theorem}

We also state, without proof, the corresponding $(p,q)$-theorem. The original $(p,q)$-theorem of Alon Kleitman \cite{AlKl} (the answer to a question of Hadwiger and Debrunner \cite{HaDe}) is about a family $\F$ of convex bodies in $\Re^d$ satisfying the $(p,q)$ property, that is, among any $p$ element of $\F$ there are $q$ that intersect. The result is that, given integers $p\ge q\ge d+1\ge 2$, there is an integer $N=N(p,q,d)$ such that for any family $\F$ satisfying the $(p,q)$ property there is a set $S\subset \Re^d$ with $|S| \le N$ such that $S\cap K \ne \emptyset$ for every $K \in \F$. In the no-dimension version the $(p,q)$ property is replaced by the $(p,q)^*$ property: among  any $p$ element of $\F$ there are $q$ that have a point in common lying in the unit ball $B$ of $\Re^d$.

\begin{theorem}\label{th:p-q} Given integers  $p\ge q\ge k \ge 2$, there is an integer $M=M(p,q,k)$ such that for any family $\F$ of convex bodies in $\Re^d$ satisfying the $(p,q)^*$ property there is a set $S\subset \Re^d$ with $|S| \le M$ such that $d(S,K) \le \frac 7{\sqrt k}$ for every $K \in \F$.
\end{theorem}

The main point here is that the bounds on $M$ and on $d(S,K)$ do not depend on $d$. Of course this is interesting only if $k>49$ as otherwise the origin is at distance 1 from every set in $\F$ intersecting $B$ (and there are at most $p-q$ sets in $\F$ disjoint from $B$). 

The rest of paper is organized the following way. Theorem~\ref{th:CalCar} is proved in Section~\ref{sec:thm1}. The next section contains the algorithmic proof of Theorem~\ref{th:CalCar2} which is another proof of Theorem~\ref{th:CalCar} with a slightly weaker constant. Section~\ref{sec:Colhel} is devoted to the proof of Theorems~\ref{th:genhelly} and~\ref{th:Colfrachelly}. The no-dimension coloured Tverberg theorem is proved in Section~\ref{sec:tverb}. Then come the proofs of the Selection Lemma and the weak $\eps$-net theorem. The last section is about Theorem~\ref{th:hellyprec} and the geometric inequality of Theorem~\ref{th:simplex}.

\section{Proof of Theorem~\ref{th:CalCar}}
\label{sec:thm1}

Given a finite set $Q \subseteq \Re^d$ denote by $\centroid(Q)$ the centroid of $Q$, that is, $\centroid(Q)=\frac 1{|Q|}\sum_{x \in Q}x$. First we prove the theorem in a special case, namely, when $a=\centroid(P_i)$ for every $i\in [r]$. Set $n_i=|P_i|$. One piece of notation: the scalar product of vectors $x,y \in \Re^d$ is written as $xy$.

We can assume (after a translation if necessary) that $a=o$. We compute the average $\Ave \centroid(Q)^2$ of $\centroid(Q)^2=\left(\frac 1{r}\sum_{x \in Q}x\right)^2=\frac1{r^2}\left(\sum_{x \in Q}x\right)^2$ taken over all transversals of the system $P_1,\ldots,P_r$. Here $(\sum_{x \in Q}x)^2$ is a linear combination of terms of the form $x^2$ ($x \in P_i),\; i\in[r]$ and $2xy$ ($x\in P_i,\,y \in P_j,\,i,j \in [r],\; i<j$). Because of symmetry, in $\Ave\,(\sum_{x \in Q}x)^2$ the coefficient of each $x^2$ with $x \in P_i$ is the same and is equal to $1/n_i$. Similarly, the coefficient of each $2xy$ with $x\in P_i,\;y\in P_j$ is the same and is equal to $1/(n_in_j)$. This follows from the fact that in every $(\sum_{x \in Q}x)^2$ out of all $x \in P_i$ exactly one $x^2$ appears, and out of all pairs $x \in P_i,y\in P_j$ exactly one $2xy$ appears. So we have
\begin{eqnarray*}
\Ave \centroid(Q)^2&=&\Ave \left(\frac 1{r}\sum_{x \in Q}x\right)^2= \frac 1{r^2}\left[\sum_1^r\frac 1{n_i}\sum_{x \in P_i}x^2+\sum_{1\le i<j\le r}\sum_{x\in P_i,\;y\in P_j}\frac 1{n_in_j}2xy\right]\\
&=&\frac 1{r^2}\left[\sum_1^r\frac 1{n_i}\sum_{x \in P_i}x^2+2\sum_{1\le i<j\le r}\left(\frac 1{n_i}\sum_{x\in P_i}x\right)\left(\frac 1{n_j}\sum_{y\in P_j}y\right)\right]\\
&=&\frac 1{r^2}\left[\sum_1^r\frac 1{n_i}\sum_{x \in P_i}x^2
+2\sum_{1\le i<j \le r}\centroid(P_i)\centroid(P_j)\right]\\
&=&\frac 1{r^2}\sum_1^r\frac 1{n_i}\sum_{x \in P_i}x^2 \le \frac 1{r^2}\sum_1^r D^2=\frac {D^2}{r},
\end{eqnarray*}
which is slightly weaker than our target. We need a simple (and probably well known) lemma.

\begin{lemma}\label{l:simple} Assume $X\subset \Re^n$, $|X|=n$, $\sum_{x\in X}x=o$ and $\diam X=D$. Then
\[
\frac 1n \sum_{x\in X}x^2 \le \frac {n-1}{2n}D^2<\frac {D^2}2.
\]
\end{lemma}

{\bf Proof.} For distinct $x,y \in X$ we have $x^2+y^2-2xy \le D^2$ and $\sum_{x\in X}x=o$ implies that $\sum_{x\in X}x^2=-\sum xy$ with the last sum taken over all distinct $x,y \in X$. Thus
\begin{eqnarray*}
\sum_{x\in X}x^2&=&-\sum xy \le \sum \frac 12 (D^2-x^2-y^2)\\
          &=& {n \choose 2}D^2-(n-1)\sum_{x\in X}x^2,
\end{eqnarray*}
implying the statement.
\qed

Using this for estimating $\Ave \centroid(Q)^2$ we get
\[
\Ave \centroid(Q)^2=\frac 1{r^2}\sum_1^r\frac 1{n_i}\sum_{x \in P_i}x^2 < \frac 1{r^2}\sum_1^r \frac {D^2}2=\frac {D^2}{2r}.
\]
This shows that there is a transversal $T$ with $|\!\centroid(T)| < D/\sqrt {2r}$. Then $d(o,\conv T) < D/\sqrt {2r}$ which proves the theorem in the special case when each $\centroid(P_i)=a$.

In the general case $a$ is a  convex combination of the elements in $P_i$ for every $i \in [r]$, that is,
\begin{equation}\label{eq:convcomb}
a=\sum_{x\in P_i}\al_i(x)x \mbox{ with }\al_i(x)\ge 0 \mbox{ and }\sum_{x\in P_i}\al_i(x)=1.
\end{equation}
By continuity it suffices to prove the statement when all $\al_i(x)$ are rational. Assume that $\al_i(x)=\frac {m_i(x)}{m_i}$ where $m_i(x)$ is a non-negative integer and $m_i=\sum_{x\in P_i}m_i(x)>0$.

Now let $P_i^*$ be the multiset containing $m_i(x)$ copies of every $x \in P_i$. Again $D=\max \diam P_i^*$, and $\centroid(P_i^*)=a$. The previous argument applies then and gives a transversal $T^*=\{p_1,\ldots,p_r\}$ of the system $P_1^*,\ldots,P_r^*$  such that
\[
d(a,\conv T^*) < \frac D{\sqrt {2r}}.
\]
To complete the proof we note that $T=T^*$ is a transversal of the system $P_1,\ldots,P_r$ as well. \qed

{\bf Remark.} One can express this proof in the following way. Choose the point $x_i \in P_i$ randomly, independently, with probability $\al_i(x)$ for all $i\in [r]$ where $\al_i(x)$ comes from (\ref{eq:convcomb}). This gives the transversal $\{x_1,\ldots,x_r\}$. We set again $a=o$. The expectation of $\left(\frac 1r (x_1+\ldots+x_r)\right)^2$ turns out to be at most
\[
\frac {D^2}{2r}\left(1-\frac 1r\sum\frac 1{n_i}\right) < \frac {D^2}{2r}.
\]
The computations are similar and this proof may be somewhat simpler than the original one. But the original one is developed further in the proof of Theorem~\ref{th:nodimColTv}. Actually, the probabilistic parts of the proofs in \cite{pisi}, \cite{Carl}, \cite{CaPa}, \cite{barman} are essentially the same except that they don't use the product distribution $\prod_1^r\al_i(\cdot)$, just the $r$-fold product of $\al(\cdot)$ coming from the convex combination $a=\sum_{x \in P}\al(x)x$.

\medskip

The above proof also works when the sets $\conv P_i$ do not intersect but there is a point close to each. Recall that $B(a,\rho)$ denotes the Euclidean ball centered at $a\in \Re^d$ of radius $\rho$.

\begin{lemma}\label{le:CalCar}
Let $P_1, \ldots, P_r$ be $r\ge 2$ point sets in $\Re^d$, $D = \max_{i\in [r]} \diam P_i$ and $\eta >0$. Assume that $B(a,\eta D)\cap \conv P_i\ne \emptyset$ for every $i \in [r]$. Then there exists a transversal $T$ such that
$$d\left( a, \conv T \right) \leq  \frac{D}{\sqrt{2r}}\sqrt{1+2(r-1)\eta^2}.$$
\end{lemma}

{\bf Proof.} The above proof works up to the point where $\sum \centroid(P_i)\centroid(P_j)$ appears. This time the sum is not zero but every term is at most $\eta^2 D^2$, and there are ${r \choose 2}$ terms. This gives the required bound. \qed

We close this section by giving a deterministic algorithm, derived by derandomizing the proof of Theorem~\ref{th:CalCar}. We state
it for the case assuming that $\alpha_i(x) = \frac{1}{n_i}$ for each $i \in [r]$ and $x \in P_i$; this is the case when $a = \centroid(P_i)$ for all $i \in [r]$. The general case follows in the same way, by derandomizing the probabilistic proof that picks each $x \in P_i$ with probability $\alpha_i(x)$ (as outlined in the equivalent formulation above).

We  will iteratively
choose the points in the sets. Assume we have selected the points
$f_i \in P_i$ for $i = s+1, \ldots, r$. We also need to be able to evaluate
the conditional expectation
$\Ex \left[ \centroid \left( \left\{ x_1, \ldots, x_s, f_{s+1}, \ldots, f_r \right\} \right)^2 \big| f_{s+1}, \ldots, f_r \right]$,
where the expectation is
over the points $x_1, \ldots, x_s$  chosen uniformly from the sets $P_1, \ldots, P_s$.
This can be done, as

\begin{eqnarray*}
\Ave &\!&\mkern-36mu \left(\frac 1{r} \left(  \sum_{i = 1}^s x_i + \sum_{i=s+1}^r f_i \right) \right)^2=\\
&=& \frac 1{r^2}\left[ \Ave \left(  \sum_{i = 1}^s x_i  \right)^2  + \left( \sum_{i=s+1}^r f_i  \right)^2 + 2 \left( \sum_{i=s+1}^r f_i  \right) \Ave \left( \sum_{i = 1}^s x_i  \right) \right].
\end{eqnarray*}
Now, as shown earlier, we have
\begin{equation} \label{eq:prefix}
\Ave \left(  \sum_{i = 1}^s x_i  \right)^2 = \sum_{i=1}^s \frac{1}{n_i} \sum_{x \in P_i} x^2.
\end{equation}
Similarly, one can compute $\Ave \left(  \sum_{i = 1}^s x_i  \right)$ exactly. Thus one can compute
$$\Ex \left[ \centroid \left( \left\{ x_1, \ldots, x_s, f_{s+1}, \ldots, f_r \right\} \right)^2 \big| f_{s+1}, \ldots, f_r \right]$$
exactly. We can pre-compute the postfix sums in equality~(\ref{eq:prefix}) at the beginning
of the algorithm, in total time $O\left( d \sum_{i}^r n_i \right)$.
Then the above expectation can be computed  in $O(1)$ time.
Now, given the sets $P_1, \ldots, P_r$, one can try all possible points  $x \in P_r$
to find the point $f_r \in P_r$ such that
$$\Ex \left[ \centroid \left( \left\{ x_1, \ldots, x_{r-1}, f_{r} \right\} \right)^2  \right] \leq \Ex \left[ \centroid \left( \left\{ x_1, \ldots, x_{r-1}, x_r \right\} \right)^2  \right], $$
in time $O(n_r)$.
This fixes the point $f_r$, and we now re-iterate to find the point $f_{r-1}$, and so on till we  have
fixed all the points $f_1, \ldots, f_r$ with the required upper-bound on $\centroid(\cdot)^2$:
$$\Ex \left[ \centroid \left( \left\{ f_1, \ldots, f_{r} \right\} \right)^2  \right]
\leq \Ex \left[ \centroid \left( \left\{ x_1, f_2, \ldots, f_r \right\} \right)^2  \right]
\leq \cdots \leq
\Ex \left[ \centroid \left( \left\{ x_1, x_2, \ldots, x_r \right\} \right)^2  \right]
< \frac{D^2}{2r}.
$$
Overall, the running time is $O\left( d \sum_{i=1}^r |P_i|\right)$.

\section{Proof of Theorem~\ref{th:CalCar2}}
\label{sec:colourfulcaratheodory}

The proof and its calculations are similar to other applications of the Frank-Wolfe method (see~\cite{bhpr, clark}
and the references therein). 
For completeness, we present the proof in our setting.

{\bf Proof.} By translation, we can assume that $a=o$.
For simpler notation we write $|q|=d(o,q)$ when $q \in \Re^d$, and the scalar product of vectors $u,v \in \Re^d$ is written as $uv$.

Initially, pick an arbitrary point of $P_1$, say $p_1 \in P_1$.
 We are going to construct a sequence $i_1=1,i_2,\ldots,i_{r-t+1}$ consisting of $r-t+1$ distinct integers with $i_j<j+t$ and a point $p_{i_j}\in P_{i_j}$ for each $i_j$ as follows. We start with an arbitrary $p_1\in P_1$ and set $q_1=p_1$. Assume $i_1=1,i_2,\ldots,i_j$ have been constructed and set $P^j=\{p_{i_1},\ldots,p_{i_j}\}$ and $I^j=\{i_1,\ldots,i_j\}$. Let $q_j$ be the nearest point of $\conv P^j$ to the origin, define $v_j=o-q_j$ and let  $Q^j$ be the union of all $P_i$ with $i \in [j+t]\setminus I^j$. Define
\[ 
p=\arg\max _{x \in Q^j} \left(x - q_j\right) v_j.
\]
Of course this point $p$ belongs to some $P_i$ with $i \in  [j+t]\setminus I^j$; 
denote it by $P_{i_{j+1}}$ and set  $p_{i_{j+1}}=p$. 
Let $q_{j+1}$ be the nearest point of $\conv P^j \cup \left\{ p_{i_{j+1}} \right\}$ to the origin.
\begin{lemma}
$$ |q_{j+1}|  \leq  \left( 1 - \frac{|q_{j}| ^2}{2 D^2} \right) \cdot |q_{j}|.$$
\label{lemma:frankwoolfonestep}
\end{lemma}

\begin{figure}
\centering
\includegraphics [scale = 0.6]{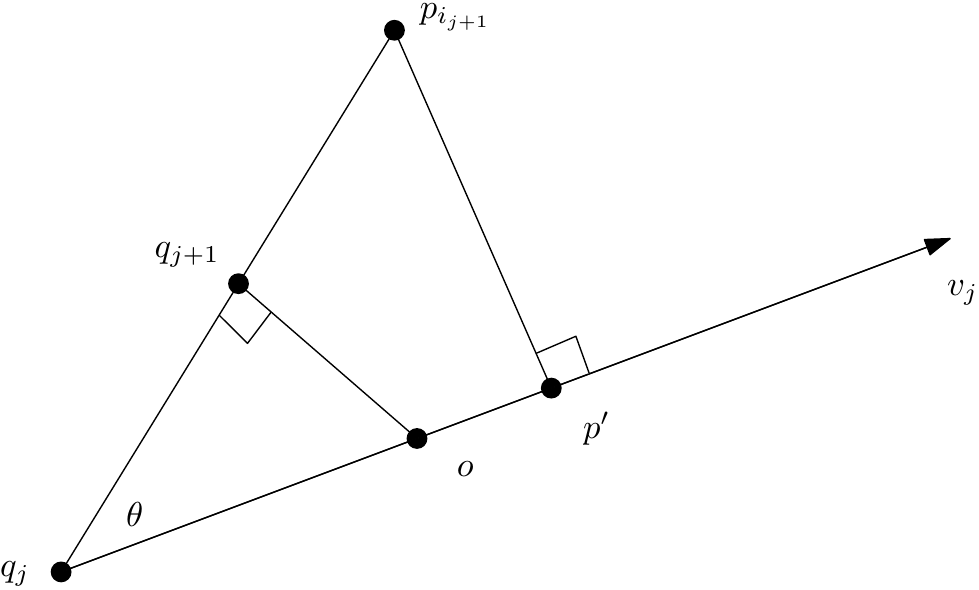}
\caption{An iteration of the algorithm.}
\label{fig:FrankWolfe}
\end{figure}

{\bf Proof.} See Figure~\ref{fig:FrankWolfe}.
In the triangle with vertices $q_j, o, q_{j+1}$, we have
$\sin \theta = \frac{|q_{j+1}|}{|q_j|}$.
On the other hand, in the triangle with vertices $q_j, p_{i_{j+1}}, p'$, we get
\begin{align*}
\cos \theta = \frac{|q_j- p' |}{|q_j -p_{i_{j+1}}|} \geq \frac{|q_j|}D.
\end{align*}
This requires $|q_j-p'| \geq |q_j|$.
Note that $o \in \conv Q^j$ due to the fact that $Q^j$ is the union of at least $t$ of the $P_i$'s. 
Thus  $|q_j-p'| \geq |q_j|$ since
any half-space containing $o$ must contain some point of $Q^j$.
Using the fact that $\sin^2 \theta + \cos^2 \theta = 1$, we get 
\begin{equation} 
|q_{j+1}| \leq \sqrt{ 1 - \left(\frac{|q_j|}{D}\right)^2} \cdot  |q_j| \leq \left(1 - \frac{|q_j|^2}{2D^2} \right) \cdot |q_j|
\end{equation}
\qed

Lemma~\ref{lemma:frankwoolfonestep} implies that at the start, when $|q_j|$ is large,
the decrease is correspondingly larger, and this slows down
with more iterations. Specifically,
for an integer $l \geq 1$, let $k_l$ be the number of indices $i_j$ such that 
\begin{equation}
\frac{D}{2^{l}} < |q_j| \leq \frac D{2^{l-1}}.
\label{eq:halving} 
\end{equation} 
Let $j' \in \left\{1, \ldots, {r-t+1} \right\}$ be the smallest index for which \eqref{eq:halving} is true.
By Lemma~\ref{lemma:frankwoolfonestep} and the fact that
$|q_j|$ is a non-increasing function of $j$, we have
\begin{align*}
|q_{j'+k_l-1}| &\leq \prod_{m=j'}^{j'+k_l-1} \left( 1 - \frac{|q_{m}|^2}{2 D^2} \right) \cdot |q_{j'}| \\
& \leq \left( e^{- \frac{|q_{j'+k_l-1}|^2}{2D^2} } \right)^{k_l} \cdot |q_{j'}|.
\end{align*}
Now the fact that $|q_{j'+k_l-1}| > \frac{D}{2^l}$ implies
that $k_l < 2^{2l} \cdot \ln 4$.
Thus the maximum number of iterations $i_j$ for which we have $|q_j| > \frac D{2^{t}}$ is
at most
\begin{align*}
\sum_{l=1}^{t} 2^{2l} \cdot \ln 4 = \ln 4 \cdot \frac{4(4^{t}-1)}{3}.
\end{align*}
This is at most $r-t+1$ for $2^t = \sqrt{ \frac{3(r-t+1)}{4 \ln 4} + 1}$.
In other words, after $r-t+1$ iterations, we have
$$|q_{r-t+1}| \leq \frac {D}{\sqrt{ \frac{3(r-t+1)}{4 \ln 4} + 1}} \leq 2 \sqrt{\frac{\ln 4}{3}} \cdot \frac D{\sqrt{r-t+1}}.$$ 
We extend this partial transversal to a complete one with an arbitrary choice of $p_i \in P_i$ for all $i\in [r-t]\setminus I^{r-t+1}$. The new transversal satisfies the required inequality.
\qed

{\bf Remark.} In this proof the first point $p_1 \in P_1$ can be chosen arbitrarily, even the condition $a \in \conv P_1$ is not needed. This implies that there are at least $|P_1|$ suitable transversals because the starting point $p_1$ can be chosen in $|P_1|$ different ways, and each gives a different transversal. We remark further that the proof is an effective algorithm that finds the transversal $Q$.

\section{Proof of Theorem~\ref{th:genhelly} and~\ref{th:frachelly}}
\label{sec:Colhel}

The {\bf proof} of Theorem~\ref{th:genhelly} goes by induction on $k$ and the case $k=1$ is trivial. For the induction step $k-1 \to k$ fix a point $q \in \Re^d$ and consider the system $\F_1,\ldots,\F_{k-1}$. By the induction hypothesis it has at least $\prod_1^{k-1}m_i$ transversal sets $\S=\{K_1,\ldots,K_{k-1}\}$ with $d(q,K(\S))>\sqrt{\rho_1^2+\ldots +\rho_{k-1}^2}$. If $K(\S)=\emptyset$, then one can extend $\S$ by any $K_k \in \F_k$ to the transversal $\T=\S \cup K_k$ that satisfies $d(q,K(\T))=\infty>\rho$. This means that $\S$ with $K(\S)=\emptyset$ can be extended to a suitable $\T$ in $n_k$ different ways.

Suppose now that $K(\S) \ne \emptyset$ and let $p$ be the point in $K(\S)$ nearest to $q$. Note that $K(\S)$ is contained in the halfspace
\[
H=\{x \in \Re^d: (p-q)(x-p)\ge 0\}.
\]
By the assumption there are at least $m_k$ sets $K\in \F_k$ with $d(p,K)> \rho_k$. For all such $K=K_k$ consider the transversal $\T=\S \cup K_k$. Then $d(q,K(\T))=\infty>\rho$ if $K(\T)=\emptyset$. Otherwise let $p'$ be the point in $K(\T)$ nearest to $q$. So $p' \in H$ and then
\[
(p'-q)^2=(p'-p)^2+(p-q)^2+2(p'-p)(p-q)> \rho_k^2+(\rho_1^2+\ldots +\rho_{k-1}^2)=\rho^2.
\]
Thus $\S$ with $K(\S)\ne\emptyset$ extends to a suitable $\T$ in $m_k$ different ways.
In both cases $\S$ can be extended to $\T$ in at least $\min \{m_k,n_k\}=m_k$ ways meaning that there are at least $\prod_1^km_i$ transversals with $d(q,K(\T))>\rho$. \qed

{\bf Proof} of Theorem~\ref{th:frachelly}. We have now a single family $\F$ with $|\F|=n$, and we assume that for every $p\in \Re^d$ there are at least $m$ sets $K\in \F$ with $d(p,K)>\rho$. We define $m=(1-\ga)n$. We want to show that $\ga \ge \be$. Fix $q\in \Re^d$ and call an ordered $j$-tuple $(K_1,\ldots,K_j)$ {\sl good} if $\bigcap_1^jK_i$ is disjoint from $B(q,\sqrt j \rho)$. We show by induction on $k$ that the number of good $k$-tuples $(K_1,\ldots,K_k)$ of $\F$ is at least
\[
(1-\ga)^k n(n-1)\ldots(n-k+1).
\]
Note that $n(n-1)\ldots(n-k+1)$ is the total number of ordered $k$-tuples of $\F$.

In the induction step of the previous proof, when considering the good $(k-1)$-tuple $(K_1,\ldots,K_{k-1})$ of $\F$ we had to consider two cases.

{\bf Case 1} when $\bigcap_1^{k-1}K_i=\emptyset$. Then we can add any $K\in \F$ distinct from $K_1,\ldots,K_{k-1}$. This is altogether $n-k+1$ good $k$-tuples of $\F$ extending the previous good $(k-1)$-tuple.

{\bf Case 2} when $\bigcap_1^{k-1}K_i\ne \emptyset$. Let $p$ be the nearest point in  $\bigcap_1^{k-1}K_i$ to $q$. Then there are $m$ sets $K\in \F$ with $d(p,K)> \rho$ and such a $K$ is different from every $K_i$ because all $K_i$ contains $p$. By the induction hypothesis $|q-p|>\sqrt {k-1}\rho$. This gives altogether $m$ good $k$-tuples of $\F$ extending the previous good $(k-1)$-tuple.

Thus each good $(k-1)$-tuple is extended to a good $k$-tuple in either $m$ or $n-k+1$ ways, finishing the induction. So the fraction of good $k$-tuples (among all $k$-tuples) is at least $(1-\ga)^k\le 1-\al$ implying that
\[
\ga \ge 1-(1-\al)^{1/k}=\be,
\]
indeed. \qed

\section{Proof of Theorem~\ref{th:nodimColTv}}\label{sec:tverb}

Before the proof of Theorem~\ref{th:nodimColTv} we need a lemma. Recall that $P$ is the disjoint union of sets (considered colours) $C_1,\ldots,C_r \subset \Re^d$, and each $C_j$ has size $k\ge 2$, the case $k=1$ is trivial.

\begin{lemma}\label{lemma:coloured} Under the above conditions there is a subset $Q\subset P$ with $|Q\cap C_j|=\lfloor \frac k2 \rfloor$ for every $j \in [r]$ such that
\begin{enumerate}
\item[(i)] $d(\centroid(Q),\centroid(P))\le \sqrt{\frac k{2(k-1)n}} \diam P$ if $k$ is even,
\item[(ii)] $d(\centroid(Q),\centroid(P))\le \sqrt{\frac {(k-2)(k+1)}{2(k-1)^2n}}\diam P$ if $k$ is odd.
\end{enumerate}
\end{lemma}

{\bf Proof.} Assume again that $\centroid(P)=o$ and write $D=\diam P$. We use an averaging argument again, this time averaging over all subsets $Q$ of $P$ with  $|Q\cap C_j|=\lfloor \frac k2 \rfloor$ for every $j \in [r]$.
We start with the case when $k$ is even.
\[
\Ave \centroid(Q)^2=\frac 4{n^2}\Ave \left(\sum_{x\in Q} x\right)^2.
\]
This is again a linear combination of terms $x^2$ for $x \in P$ and $2xy$ for $x,y\in P$, $x\ne y$ of distinct colours and  $2xy$ for $x,y\in P$, $x\ne y$ of the same colour. It is clear that each $x^2$ goes with coefficient $\frac 12$, each $2xy$ from different colours with coefficient $\frac 14$ while the coefficient of $2xy$ with $x,\, y$ of the same colour (and $x\ne y$) is
\[ \frac {{k/2 \choose 2}}{{k \choose 2}} =\frac {k-2}{4(k-1)}=\frac 14\left(1-\frac 1{k-1}\right).
\]
Thus, writing $\sum\,\! _{(1)}$ resp. $\sum\,\! _{(2)}$ for the sum taken over pairs $x,y$ of distinct colour and of the same colour,
\begin{eqnarray*}
\Ave \left(\sum_{x\in Q} x\right)^2&=&\frac 12 \sum_{x \in P}x^2+ \frac 14  \sum\,\! _{(1)} 2xy +\frac 14 \left(1-\frac 1{k-1}\right)\sum\,\! _{(2)} 2xy\\
&=&\frac 14 \sum_{x \in P}x^2+\frac 14 \left(\sum_{x\in P}x\right)^2-\frac 1{4(k-1)}\sum\,\! _{(2)} 2xy\\
&=&\frac 14\left(1-\frac 1{k-1}\right) \sum_{x \in P}x^2 - \frac 1{4(k-1)}\sum_{j=1}^k\left(\sum_{x \in C_j}x\right)^2\\
&\le& \frac 14 \frac k{k-1}  \sum_{x \in P}x^2\le  \frac 14 \frac {k}{k-1} \frac {nD^2}{2}.
\end{eqnarray*}
according to Lemma~\ref{l:simple}. Returning now to $\Ave \centroid(Q)^2$ we have
\[
\Ave \centroid(Q)^2=\frac 4{n^2}\Ave \left(\sum_{x\in Q} x\right)^2 \le \frac {kD^2}{2(k-1)n}=\frac 1{k-1}\frac {D^2}{2r}.
\]
Consequently there is a $Q \subset P$ satisfying $(i)$.
Assume next that $k$ is odd: $k=2s+1$, say, with $s \ge 1$. So the average is to be taken over all subsets $Q$ of $P$ with $|Q\cap C_j|=s$ for every $j \in [k]$, and
\[
\Ave \centroid(Q)^2=\frac 1{(sk)^2} \Ave \left(\sum_{x\in Q} x\right)^2.
\]
The coefficients of $x^2$ and $2xy$ in $\Ave \left(\sum_{x\in Q} x\right)^2$ are determined the same way as before and we have
\begin{eqnarray*}
\Ave \left(\sum_{x\in Q} x\right)^2&=&\frac s{2s+1}\sum_{x\in P}x^2+\frac {s^2}{(2s+1)^2} \sum\,\!_{(1)} 2xy +\frac {s(s-1)}{(2s+1)2s}\sum\,\!_{(2)} 2xy\\
&=& \frac s{2s+1}\left(\sum_{x\in P}x^2+\frac s{2s+1}\sum\,\!_{(1)}2xy+\frac{s-1}{2s}\sum\,\!_{(2)} 2xy\right)\\
&=& \frac s{2s+1}\left[\frac {s+1}{2s+1}\sum_{x\in P}x^2+\frac s{2s+1}\left(\sum_{x\in P}x\right)^2+\left(\frac {s-1}{2s}-\frac s{2s+1}\right)\sum\,\!_{(2)} 2xy \right]\\
&=& \frac {s(s+1)}{(2s+1)^2}\left[\sum_{x\in P}x^2-\frac 1{2s}\sum\,\!_{(2)} 2xy\right]\\
&=& \frac {s(s+1)}{(2s+1)^2}\left[ \left(1-\frac 1{2s}\right)\sum_{x\in P}x^2-\frac 1{2s}\sum_1^k\left(\sum_{x \in C_j}x\right)^2 \right]\\
&\le& \frac {s(s+1)}{(2s+1)^2}\left(1-\frac 1{2s}\right)\sum_{x\in P}x^2\le
 \frac {(s+1)(2s-1)}{2(2s+1)^2}\frac {nD^2}2,
\end{eqnarray*}
where Lemma~\ref{l:simple} was used again. Here $n=rk=r(2s+1)$ so
\[
\Ave \centroid(Q)^2 \le \frac 1{s^2k^2}\frac {(s+1)(2s-1)}{2(2s+1)^2}\frac {nD^2}2=\frac {(k+1)(k-2)}{k(k-1)^2} \frac {D^2}{2r} < \frac 1{k-1}\frac {D^2}{2r},
\]
where the last inequality follows easily from $k\ge 3$. This proves part $(ii)$.\qed

\begin{corollary}\label{cor:ColPartition} Under the conditions of Lemma~\ref{lemma:coloured} there is a partition $Q_0,Q_1$ of $P$ with  $|Q_0\cap C_j|=\lfloor \frac k2 \rfloor $ and $|Q_1\cap C_j|=\lceil \frac k2 \rceil$ for every $j \in [r]$ such that $|\!\centroid(Q_0)|=|\!\centroid(Q_1)|$ when $k$ is even, and $\frac {k-1}{k+1}|\!\centroid(Q_0)|=|\!\centroid(Q_1)|$ when $k$ is odd. Moreover
\[
d(\centroid(Q_1),\centroid(P))\le d(\centroid(Q_0),\centroid(P)) \le \frac 1{\sqrt{k-1}}\frac D{\sqrt {2r}}.
\].
\end{corollary}

{\bf Proof.} Set $Q_0=Q$ where $Q$ comes from Lemma~\ref{lemma:coloured} and $Q_1=P\setminus Q_0$. In the even case $\centroid(Q_0)+\centroid(Q_1)=2\centroid(P)=o$ again. In the odd case $s\centroid(Q_0)+(s+1)\centroid(Q_1)=\centroid(P)=o$, implying that $|\!\centroid(Q_1)|=\frac {k-1}{k+1}|\!\centroid(Q_0)|$. Moreover, $n=rk$ and in all cases $|\!\centroid(Q_1)|\le |\!\centroid(Q_0)|\le \frac 1{\sqrt{k-1}}\frac D{\sqrt {2r}}$.  \qed

{\bf Proof} of Theorem~\ref{th:nodimColTv}. We build an incomplete binary tree. Its root is $P$ and its vertices are subsets of $P$. The children of $P$ are $Q_0,Q_1$ from the above Corollary, the children of $Q_0$ resp. $Q_1$ are $Q_{00},Q_{01}$ and $Q_{10},Q_{11}$ obtained again by applying Corollary~\ref{cor:ColPartition} to $Q_0$ and $Q_1$. We split the resulting sets into two parts of as equal sizes as possible the same way, and repeat. We stop when the set $Q_{\de_1\ldots\de_h}$ contains exactly one element from each colour class. In the end we have sets $P_1,\ldots,P_r$ at the leaves. They form a partition of $P$ with $|P_i\cap C_j|=1$ for every $i \in [r]$ and  $j \in [k]$. We have to estimate $d(\centroid(P_i),\centroid(P))$. Let $P,Q^1,\ldots,Q^h,P_i$ be the sets in the tree on the path from the root to $P_i$. Using the Corollary gives
\begin{eqnarray*}
d(\centroid(P),\centroid(P_i))&\le& d(\centroid(P),\centroid(Q^1))+d(\centroid(Q^1),\centroid(Q^2))+\ldots+d(\centroid(Q^{h}),\centroid(P_i))\\
 &\le& \left[ \frac 1{\sqrt{k-1}}+ \frac 1{\lfloor \sqrt {k/2 \rfloor-1}} + \frac 1{\sqrt {\lfloor k/4 \rfloor-1}}+\ldots \right]\frac D{\sqrt {2r}}\\
 &\le&  (1+\sqrt 2) \frac D{\sqrt {r}},
\end{eqnarray*}
as one can check easily. \qed

We mention that with a little extra care the constant $1+\sqrt 2=2.4142..$ can be brought down to 2.02.

\section{Proofs of the No-Dimension Selection and  Weak $\eps$-net Theorems}

{\bf Proof} of Theorem~\ref{th:selection}. This is a combination of Lemma~\ref{le:CalCar} and the no-dimension Tverberg theorem, like in \cite{bar}. We assume that  $n=kr+s$ with $0\le s\le r-1$ ($k$ an integer) and set $\ga=2+\sqrt 2$. The no-dimension Tverberg theorem implies that $P$ has a partition $\{P_1,\ldots,P_k\}$ such that $\conv P_i$ intersects the ball $B\left(q,\ga \frac {D}{\sqrt r}\right)$ for every $i \in [k]$ where $q \in \Re^d$ is a suitable point.

Next choose a sequence $1\le j_1\le j_2 \le \ldots \le j_r \le k$ (repetitions allowed) and apply Lemma ~\ref{le:CalCar} to the sets $P_{j_1},\ldots,P_{j_r}$, where we have to set $\eta= \frac{\ga}{\sqrt r}$. This gives a transversal $T_{j_1\ldots j_r}$ of $P_{j_1},\ldots,P_{j_r}$ whose convex hull intersects the ball
\[
B\left(q,\frac {D}{\sqrt{2r}}\sqrt{1+2(r-1) \eta^2} \right).
\]
The radius of this ball is
\[
\frac D{\sqrt {2r}} \sqrt{1 +2\frac {r-1}{r}\ga^2} \le \frac D{\sqrt r}\sqrt{6.5+4\sqrt 2}<3.5\frac D{\sqrt r}
\]
as the function under the first large square root sign is decreasing with $r$ and for $r=2$ it is $13+8\sqrt 2$. So the convex hull of all of these transversals intersects $B\left(q,\frac {3.5D}{\sqrt r}\right)$. They are all distinct $r$-element subsets of $P$ and their number is
\[
{k+r-1 \choose r}={\frac {n-s}r +r-1 \choose r} \ge r^{-r}{n \choose r},
\]
as one can check easily. \qed

The {\bf proof} of Theorem~\ref{th:epsnet} is an algorithm that goes along the same lines as in the original weak $\eps$-net theorem \cite{ABFK}. Set $F:=\emptyset$ and let $\H$ be the family of all $r$-tuples of $P$. On each iteration we will add a point to $F$ and remove $r$-tuples from $\H$.

If there is $Y \subset P$ with $\left(F+\frac {3.5D}{\sqrt r}B\right)\cap \conv Y = \emptyset$, then apply Theorem~\ref{th:selection} to that $Y$ resulting in a point $q\in \Re^d$ such that the convex hull of at least
\[
\frac 1{r^r}{\eps n \choose r}
\]
$r$-tuples from $Y$ intersect $B\left(q,\frac {3.5D}{\sqrt r}\right)$. Add the point $q$ to $F$ and delete all $r$-tuples $Q\subset Y$ from $\H$ whose convex hull intersects $B\left(q,\frac {3.5D}{\sqrt r}\right)$. On each iteration the size of $F$ increases by one, and at least $r^{-r}{\eps n \choose r}$ $r$-tuples are deleted from $\H$. So after
\[
\frac {{n \choose r}}{\frac 1{r^r}{\eps n \choose r}} \le \frac {r^r}{\eps^r}
\]
iterations the algorithm terminates as there can't be any further $Y \subset P$ of size $\eps n$ with  $\left(F+\frac{3.5D}{\sqrt r}B\right)\cap \conv Y = \emptyset$. Consequently the size of $F$ is at most $r^r{\eps^{-r}}$.\qed

\section{Proof of Theorem~\ref{th:hellyprec}}\label{sec:simplex}

We begin with a simple observation.
\begin{proposition}\label{prop:boundary}
Consider a finite family $\F$ of convex bodies in $\Re^d$, a point $p$ in $\Re^d$ and a natural number $k$ at most $d+1$. Assume that every point is at positive distance from at least one of the elements of $\F$. If $\F^*$ is a subfamily of size $k$ such that the distance of its intersection from $p$ is maximal among all families of size $k$, then the closest point $q$ to $p$ in the intersection lies in the intersection of the respective boundaries.
\end{proposition}

{\bf Proof.} Assume the contrary. Then the distance from $p$ to the intersection over $\F^*$ is attained at a subfamily $\F'$ of size $k-1$. By assumption, there exists an element $K$ of $\F$ that does not contain $q$. The subfamily $\F'\cup\{K\}$ has $k$ elements and its intersection is farther from $p$ than $q$.\qed

Next we prove the geometric inequality about simplices.

{\bf Proof} of Theorem~\ref{th:simplex} Let  $p_1$, $\dots$, $p_n$ be points in general position in $\Re^{n-1}$, their convex hull is a simplex $\tri$ whose inradius is $r$. For each $j=[n]$ denote by $\sigma_j$ the facet of $\tri$ opposite to $p_j$.

We proceed by induction on $n$. For $n=k+1$ the statement is tautological. Let now $h_j$ be the height of $p_j$ over $\sigma_j$ and denote by $T_j$ the $(n-2)$-dimensional volume of $\sigma_j$. Calculating the volume of $\tri$ from these heights and the inradius, respectively, we get that for each $j$ we have $r\sum_i T_i=h_jT_j$ and consequently $h_j=\frac{\sum_i T_i}{T_j}r$.

For any fixed $j$, consider the slice of the inscribed ball parallel to $\sigma_j$ at height $z\in[0,2r]$ over this facet. This slice is a ball of radius $\sqrt{r^2-(z-r)^2}$; it lies entirely in the simplex, so its stereographic projection from the vertex $p_j$ onto the hyperplane of $\sigma_j$ lies entirely in $\sigma_j$ and thus the radius of the projection is a lower bound on the inradius of $\sigma_j$. The radius of the projection is $\sqrt{r^2-(z-r)^2}\frac{h_j}{h_j-z}$; for fixed $r$ and $h_j$ the maximum is attained at $z=\frac{h_j}{h_j-r}r$ and has value $\varrho_j := r\sqrt{\frac{h_j}{h_j-2r}}=r\sqrt\frac{\sum_i T_i}{\sum_i T_i - 2T_j}$. By the induction hypothesis this implies that any ball that meets the affine span of each $(k-1)$-dimensional face of $\sigma_j$ has radius at least $\lambda_{n-1}\varrho_j$.

Assume now that a ball of radius $R$ meets the affine span of each $(k-1)$-dimensional face of $\tri$; let $m_j$ be the (signed) height of its center above $\sigma_j$. By computation of volume of the simplex we have $\sum_i m_iT_i=r\sum_i T_i$. By the induction hypothesis, in order to meet the affine span of each $k$-dimensional face of $\sigma_j$ in particular the intersection of the ball with $\sigma_j$ -- an $(n-1)$-ball of radius $\sqrt{R^2-m_j^2}$ -- has to have radius at least $\lambda_{n-1}\varrho_j$, hence $R \geq \sqrt{m_j^2+\lambda_{n-1}^2\varrho_j^2}$ holds for all $j$. Introducing the notation $\alpha_j:= \frac{T_j}{\sum_i T_i}$, we claim that the three conditions
\[
\sum_{i=1}^n \frac{m_i}{r}\alpha_i=1, \;\; \sum_{i=1}^n \alpha_i=1,\mbox{ and } \ 0\leq \alpha_i \leq \frac{1}{2} \mbox{ for all } i\in [n]
\]
imply that
\begin{equation}\label{eq:opt}
\max_j \left\{m_j^2+\lambda_{n-1}^2\varrho_j^2\right\} = r^2 \max_j \left\{\left(\frac{m_j}{r}\right)^2+\lambda_{n-1}^2\frac{1}{1-2\alpha_j} \right\} \geq \lambda_n^2r^2
\end{equation}
This will finish the proof of the induction step.

To prove the inequality \eqref{eq:opt}, form the weighted average of the expressions $\left(\frac{m_j}{r}\right)^2+\lambda_{n-1}^2\frac{1}{1-2\alpha_j}$ with weights $\alpha_j$:
$$
\sum_{j=1}^n \alpha_j \left(\left(\frac{m_j}{r}\right)^2+\lambda_{n-1}^2\frac{1}{1-2\alpha_j}\right) = \sum_{j=1}^n \alpha_j \left(\frac{m_j}{r}\right)^2 + \lambda_{n-1}^2 \sum_{j=1}^n \frac{\alpha_j}{1-2\alpha_j}.
$$
By convexity of $t\mapsto t^2$ the first sum -- considered as a weighted average -- is at least $\left(\sum_j\alpha_j\frac{m_j}{r}\right)^2=1$ and by convexity of the function $t\mapsto \frac{t}{1-2t}$ the second sum -- considered as a regular average -- is at least
$$
n\lambda_{n-1}^2\frac{\frac{\sum_j\alpha_j}{n}}{1-2\frac{\sum_j\alpha_j}{n}} = \lambda_{n-1}^2\frac{n}{n-2} = \frac{n(n-k-1)}{k} = \lambda_n^2-1.
$$
Hence the sum of the two parts is at least $\lambda_n^2$ and consequently at least one of the weighted summands is at least $\lambda_n^2$.
This proves (\ref{eq:opt}) and finishes the proof.\qed

We need a slight strengthening of this inequality. Given the simplex $\tri$, let $G_i$ denote the closed halfspace satisfying $\tri \cap \G_i=\sigma_i$. We define $C(J)$, the cone over the $(k-1)$-face $\conv\{p_j: j\in J\}$ as $C(J)=\bigcap_{j\in J} G_j$.

\begin{lemma}\label{strong} Under the conditions of Theorem~\ref{th:simplex}, any ball intersecting $C(J)$ for every $J \subset [n], |J|=k$ has radius at least $\lambda_n r$.
\end{lemma}
 The {\bf proof} follows directly from Proposition~\ref{prop:boundary}: the intersection of the boundaries of the halfspaces $G_j, j\in J$
is exactly the affine hull of the corresponding $(k-1)$-face. \qed

\medskip
{\bf Proof} of Theorem~\ref{th:hellyprec}. First we show how to replace each $K_i$ with a polytope.
Choose a point $z(J) \in K(J) \cap B(b,1)$ for every $J \subset [n], |J|=k$ and set $K_i^*=\conv\{z(J): i \in J\}$. The new family $\F^*$ satisfies the same conditions as $\F$, each $K_i^*$ is a polytope and is a subset of $K_i$. Thus if no ball $B(q,1/\sqrt k)$ intersects all $K_i \in \F$, then it does not intersect all $K_i^* \in \F^*$ either.

Next set $r=\sqrt{\frac{n-k}{k(n-1)}}$ and define $K_i=K_i^*+rB$; we have to show that $\bigcap_1^nK^i \ne \emptyset$. Assume the contrary. Then there are closed halfspaces $H_i$
in general position such that $K_i \subset H_i$ and  $\bigcap_1^nH^i = \emptyset$. Write $a_i$ for the outer (unit) normal of $H_i$ and $A$ for the linear span of $a_1,\ldots,a_n$.
It is clear that $A$ is a copy of $\Re^{n-1}$. Let $H_i^r$ denote the halfspace contained in $H_i$ such that their bounding hyperplanes are exactly a distance $r$ apart. Then
\[
\tri=A \cap \bigcup_1^nH_i^r
\]
is a simplex in $A\cong \Re^{n-1}$ whose inradius is at least $r$. The outer cone of $\tri$ over the face $\conv\{v_j:j\in J\}$ is $C(J)=A\cap \bigcap_{j \in J}H_j^r$. Lemma~\ref{strong} applies now and shows that for every $q \in A$ one of the outer cones with $|J|=k$ is farther than $\lambda_nr=1$. A contradiction with the assumption that $K(J)$ has a point in $B(b,1)$.\qed

\bigskip
{\bf Acknowledgements.}  K.\ A.\ was supported by ERC StG 716424 - CASe and ISF Grant 1050/16. I.\ B.\ was supported  by the Hungarian National Research,
Development and Innovation Office NKFIH Grants K 111827 and K 116769, and by ERC-AdG 321104. N.\ M.\ was supported by the grant ANR SAGA (JCJC-14-CE25-0016-01).
T.T. was supported by the Hungarian National Research,
Development and Innovation Office NKFIH Grants NK 112735 and K 120697.

\noindent
Karim Adiprasito\\
Einstein Institute for Mathematics,\\
Hebrew University of Jerusalem\\
Edmond J. Safra Campus, Givat Ram\\
91904 Jerusalem, Israel\\
{\tt email: adiprasito@math.huji.ac.il}

\medskip

\noindent
Imre B\'ar\'any \\
Alfr\'ed R\'enyi Institute of Mathematics,\\
Hungarian Academy of Sciences\\
13 Re\'altanoda Street Budapest 1053 Hungary\\
and\\
Department of Mathematics\\
University College London\\
Gower Street, London, WC1E 6BT, UK\\
{\tt barany.imre@renyi.mta.hu}\\

\medskip
\noindent
Nabil H. Mustafa \\
Universit\'e Paris-Est, \\
Laboratoire d'Informatique Gaspard-Monge, Equipe A3SI,\\
ESIEE Paris.\\
{\tt mustafan@esiee.fr}

\medskip
\noindent
Tam\'as Terpai\\
Department of Analysis, Lorand Eotvos University\\
P\'azm\'any P\'eter s\'et\'any 1/C, Budapest\\
H-1053 Hungary\\
{\tt terpai@math.elte.hu}

\begin{thebibliography}{99}

\bibitem{AlonLee} N. Alon, T. Lee, A. Shraibman, S. Vempala, 
The approximate rank of a matrix and its algorithmic applications: approximate rank, 
{\it 45-th Symposium on the Theory of Computing (STOC)}, 675-684, 2013.

\bibitem{ABFK} N. Alon, I. B\'ar\'any, Z. F\"uredi, and D. Kleitman, Point selections and weak $\epsilon $--nets for convex hulls, {\it Combinatorics, Probability, and Computation}, {\bf 1} (1992), 189--200.

\bibitem{AlKl} N. Alon, and D.J. Kleitman, Piercing convex sets and the Hadwiger-Debrunner $(p, q)$-problem, {\it Adv. Math.}, {\bf 96} (1992), 103--112.

\bibitem{ABBFM} L. J. Arocha, I. B{\'a}r{\'a}ny, J. Bracho, R. Fabila, and L. Montejano, Very Colorful Theorems, {\it Discrete Comput. Geom.} {\bf 42} (2009), 142--154.

\bibitem{bar} I. B\'ar\'any,  A generalization of Charath\'eodory's theorem, {\it Discrete Math.} {\bf 40} (1982), 141--152.

\bibitem{BF} I. B\'ar\'any,  Z.  F\"uredi: Computing the volume is difficult, {\it Discrete and
Comput. Geometry} {\bf 2} (1987), 319--326.

\bibitem{BFMOP}  I. B\'ar\'any, F. Fodor, L. Montejano, D. Oliveros, A. P\'or: Colourful and fractional (p,q) theorems, {\it Discrete Comp. Geom.},  {\bf 51} (2014), 628--642.

\bibitem{barman} S. Barman, Approximating Nash equilibria and dense bipartite subgraphs via an approximate version of Carath\'eodory's theorem, STOC'15—Proceedings of the 2015 ACM Symposium on Theory of Computing, 361--369, ACM, New York, 2015.

\bibitem{BlaZie} P. M. Blagojevi\'c, B. Matschke, G. M. Ziegler, Optimal bounds for the colored Tverberg problem, {\it J. Eur. Math. Soc. (JEMS)}, {\bf 17} (2015), 739--754.

\bibitem{bhpr} A. Blum, S. Har-Peled, B. Raichel, Sparse approximation via generating point sets, 
 Proceedings of the Twenty-seventh Annual ACM-SIAM Symposium on Discrete Algorithms, 548--557, 2016.


\bibitem{Carath} C. Carath\'eodory, \"Uber den Variabilit\"atsbereich der Koeffizienten von Potenzreihen,
{\it Math. Annalen}, {\bf 64} (1907), 95--115.

\bibitem{Carl} B. Carl. Inequalities of Bernstein-Jackson-type and the degree of compactness of operators
in Banach spaces, {\it Annales de l’institut Fourier}, {\bf 35} (1985), 79--118.

\bibitem{CaPa} B. Carl, A. Pajor, Gel'fand numbers of operators with values in a Hilbert space, {\it Invent. Math.} {\bf 94} (1988), 479--504.

\bibitem{Cas} J. W. S. Cassels, Measures of the non-convexity of sets and the Shapley–Folkman–Starr
theorem, {\it Math. Proc. Cambridge Philos. Soc.}, {\bf 78} (1975), 433--436.

\bibitem{clark} K. L. Clarkson, Coresets, sparse greedy approximation, and the Frank-Wolfe algorithm, {\it ACM Trans. Algorithms}, {\bf 6} (2010), 30 pp.

\bibitem{FMMZ} M. Fradelizi, M. Madiman, A. Marsiglietti, A. Zvavitch, The convexification effect of Minkowski summation, {\it EMS Surv. Math. Sci.}, {\bf 5} (2018), 1--64.

\bibitem{GiMi} A. A. Giannopoulos and V. Milman, Concentration property on probability spaces, {\it Adv. Math.} {\bf 156} (2000), 77--106.

\bibitem{Gue} O. Gu\'edon Concentration phenomena in high dimensional geometry, {\it Journées MAS 2012}, 47--60,
ESAIM Proc., {\bf 44}, EDP Sci., Les Ulis, 2014.

\bibitem{FrWo} M. Frank and Ph. Wolfe. An algorithm for quadratic programming, {\it Naval Res. Logist. Quart.}, {\bf 3} (1956), 95--110.

\bibitem{HaDe} H. Hadwiger, H. Debrunner, \"Uber eine Variante zum Hellyschen Satz, {\it Arch.
Math.}, {\bf 8} (1957), 309--313.

\bibitem{HoPaTv} A. Holmsen, J. Pach, H. Tverberg, Points surrounding the origin, {\it Combinatorica}, {\bf 28} (2008), 633--644.

\bibitem{jun} H. W. Jung, \"Uber die kleinste Kugel, die eine r\"aumliche Figur einschliesst, {\it J. Reine Angew. Math.}, {\bf 123} (1901), 241--257.

\bibitem{Kalai} G.  Kalai, Intersection  patterns  of  convex  sets, \emph{Isr. J. Math.}, \textbf{48} (1984), 161--174.

\bibitem{KaLi} M. Katchalski, A. Liu, A problem of geomtrey in $\Re^n$, \emph{Proc. AMS}, {\bf 75} (1979), 284--288.

\bibitem{Mat} J. Matou\v{s}ek, Lectures on discrete geometry, Springer (2002), New York.

\bibitem{MKim}
M. Kim, A note on the colorful fractional Helly theorem, \emph{Discrete Math.}, 340 (2017), 3167-–3170.

\bibitem{pisi} G. Pisier, Remarques sur un r\'esultat non publi\'e de B. Maurey, {\it S\'eminaire Analyse fonctionnelle}, (1981), 1–12.

\bibitem{Rad} R. Rado, A theorem on general measure, {\it J. London Math. Soc.} {\bf 41} (1966) 123--128.

\bibitem{Sta} R. M. Starr, Quasi-equilibria in markets with non-convex preferences, {\it Econometrica}, {\bf 37} (1969), 25--38. 

\bibitem{Sta2} R. M. Starr, Approximation of points of convex hull of a sum of sets by points of the sum: an elementary approach, {\it J. Econom. Theory}, {\bf 25} (1981), 314--317.

\bibitem{Tal} M. Talagrand, Concentration of measure and isoperimetric inequalities in product spaces, {\it Inst. Hautes Études Sci. Publ. Math.}, {\bf 81} (1995), 73--205.

\bibitem{tver} H. Tverberg, A generalization of Radon's theorem, {\it J. London Math. Soc.} {\bf 21} (1946) 291--300.

\bibitem{ZiVr} R. \v{Z}ivaljevi\'c, S. Vre\'cica, The colored Tverberg's problem and complexes of injective functions, {\it J. Combin. Theory Ser. A}, {\bf 61} (1992), 309--318.

\end{thebibliography}
\end{document}